\documentclass[11pt]{amsart}

\usepackage{amsmath}
\usepackage{amssymb}
\usepackage{amsfonts}
\usepackage{amsthm}
\usepackage{enumerate}
\usepackage{verbatim}
\usepackage{color}
\usepackage[mathscr]{eucal}

\usepackage{graphicx,epsfig}

\def\inn#1#2{\langle#1,#2\rangle}   
\renewcommand{\d}{\,{\rm d}}

\begin{document}

\nocite{*}

\title{Some Problems in Harmonic Analysis}
\author []{With contributions by \\ \\Almut Burchard, Ciprian Demeter, \\Michael Goldberg,  Alexander Kiselev, \\Ben Krause, Detlef  M\"uller, \\Alexander  Reznikov, Andreas Seeger, \\Christoph Thiele, James Wright}

\thanks{
The conference was funded by a  conference grant and a Research Training Group grant of  the National Science  Foundation and by  grants of  the Institute of Mathematics and Applications, the Fields Institute and the Pacific Institute for the Mathematical Sciences. The organizers thank the staff of the Department of Mathematics of the University of Wisconsin for providing valuable logistical support}

\address{L. Grafakos, Department of Mathematics\\
University of Missouri\\
Columbia, MO 65211\\ 
USA}
\email{{\tt grafakosl@missouri.edu}}

\address{D. Oliveira e Silva, Mathematisches Institut 
der Universit\"at Bonn,
Endenicher Allee 60,
53115 Bonn, Germany}
\email{{\tt dosilva@math.uni-bonn.de }}
 
\address{M. Pramanik, Department of Mathematics\\
University of British Columbia\\
Room 121, 1984 Mathematics Road\\
Vancouver, B.C.,  \\  
Canada V6T 1Z2 }
\email{{\tt malabika@math.ubc.ca}}

\address{A. Seeger, Department of Mathematics\\
University of Wisconsin-Madison\\
613 Van Vleck Hall, 480 Lincoln Drive\\
Madison, WI 53706   \\ 
USA}
\email{{\tt seeger@math.wisc.edu}}

\address{B. Stovall, Department of Mathematics\\
University of Wisconsin-Madison\\
721 Van Vleck Hall, 480 Lincoln Drive\\
Madison, WI 53706   \\ 
USA}
\email{{\tt stovall@math.wisc.edu}}

\maketitle

\bigskip

\bigskip

\bigskip

In May 2016, we organized a conference in harmonic analysis  in honor of Professor Michael Christ, 
 on the campus of the  University of Wisconsin in Madison.
 We   are happy to present   sixteen open problems, almost all of which
    were contributed by participants of a problem session held in the afternoon of May 19, 2016.

 \bigskip

\noindent   Loukas Grafakos \\
  Diogo Oliveira e Silva \\
   Malabika Pramanik \\
    Andreas Seeger \\
     Betsy Stovall 
 
\newpage
\section{A Convolution Inequality on Spheres}
\begin{center}{\it Contributed by Almut Burchard }\end{center}
\bigskip

How can one extend the  
{\it Riesz--Sobolev inequality}
\begin{equation} \label{eq:RS}
\int_{\mathbb{R}^d} \int_{\mathbb{R}^d}
f(x) g(x-y) h(y)\, \d x\d y\le
\int_{\mathbb{R}^d} \int_{\mathbb{R}^d}
f^{\#}(x) g^{\#}(x-y) h^{\#}(y)\, \d x \d y
\end{equation}
to spheres (and perhaps more general
homogeneous spaces)? Here, $f,g$ are nonnegative functions
on $\mathbb{R}^d$
that decay at infinity, and $f^{\#},g^{\#}$
are their symmetric decreasing rearrangements.
The simplest case,
where  the middle function $g(x-y)=g^{\#}(x-y)$
is a decreasing function of the
distance between $x$ and $y$, extends
directly from $\mathbb{R}^d$ to the standard spheres
and hyperbolic spaces. 
The full inequality, where all three functions may vary,
is equivalent by a convex duality argument to the statement that
\begin{equation}\label{eq:Phi}
\int \Phi(f*g) \le \int \Phi(f^{\#}*g^{\#})
\end{equation}
for every increasing convex function $\Phi$ on $\mathbb{R}_+$.
Eq.~\eqref{eq:RS} is a refinement of the Brunn--Minkowski inequality 
and its proofs use the affine geometry of $\mathbb{R}^n$.
A folklore result says that Eq.~\eqref{eq:RS} fails for
the natural convolution operator on the orthogonal group
$O(d+1)$ (or the sphere $S^d$), given by
$$ (f*g)(Q) = \int_{O(d+1)} f(y) g(Q^{-1}y),$$
where integration is with respect to the Haar measure,
and $f^{\#}$ denotes symmetric decreasing
rearrangement about the unit element (or the
north pole on the sphere).
(This result does not seem to exist in the literature.
For $d=2$, the inequality would imply a version of the isoperimetric
inequality on $SO(3)= U(2)/\{\pm\}$ which is false for subsets
whose complement is small.)

\bigskip 
Consider restricting the convolution from
the orthogonal group to the conjugacy class of reflections.
For any vector $u\in S^d$ (viewed as the unit
sphere in $\mathbb{R}^{d+1}$),
denote by $\sigma_u$ the Euclidean reflection 
on the hyperplane $u^\perp\subset \mathbb{R}^{d+1}$.  
Define a {\em convolution
operator} for functions $f,g$ on $S^d$ by
$$
(f*g)(u) = \int_{S^{d}} f(x)g(\sigma_u(x))\, \d x\,.
$$

 \medskip
\noindent {\bf Problem 1:} {\it 
Given the above convolution on the sphere, does Eq.~\eqref{eq:Phi} hold for
every convex increasing
function $\Phi$ on $\mathbb{R}_+$?
Here, $f^{\#}, g^{\#}$ are the symmetric
decreasing rearrangements of $f$ and $g$, 
i.e., $f^{\#}$, $g^{\#}$ are equimeasurable
with $f,g$, and are non-increasing functions
of the distance from the north pole.
}
\bigskip

Clearly,
$f*g$ is symmetric under the antipodal map $u\to -u$.
Note also that $f^{\#}*g^{\#}$ is
symmetric decreasing about the equator.
By duality, an equivalent question is whether
$ \int\bigl(f*g\bigr)h
$
increases if $f,g$ are
replaced by $f^\#,g^\#$,
and $h$ is rearranged to be symmetric decreasing 
about the equator.

\bigskip
The question is open in any dimension $d>1$.
On $S^1$, Eq.~\eqref{eq:RS} reduces to a theorem
of A. Baernstein~\cite{[1]}, because the
convolution operator agrees with the standard one
except for a dilation.

The convolution operator appears in connection
with two-point symmetrizations on the sphere~\cite{[4]}, \cite{[5]}.
To provide some evidence in favor of the inequality, 
consider the case $\Phi(t)=t$. 
Fubini's theorem and the change of variables $y=\sigma_u(x)$ 
yield
$$
\int_{S^d} (f*g)(u) \, \d u =
\int_{S^d}\int_{S^d} f(x)g(y) |x-y|^{-(d-1)}\,  \d x\d y\,,
$$
which increases under symmetrization, since
the Jacobian $|x-y|^{-(d-1)}$ is a symmetric decreasing function
of the distance. 

\bigskip

\bigskip

\address{A. Burchard, 
Department of Mathematics, University of Toronto, 
Toronto, Ontario, Canada M5S 2E4.}

{\it E-mail address:} 
\email{\tt almut@math.toronto.edu}

\newpage

\section{Questions about generalized additive energy}

\begin{center} {\it Contributed by Ciprian Demeter} \end{center}

\bigskip\bigskip

For sets $A\subset \mathbb R^d$ consider the generalized additive energy  
$$
{\mathbb E}_n(A)=\#\{(a_1,\ldots,a_{2n})\in A^{2n}:a_1+\ldots+a_n=a_{n+1}+\ldots +a_{2n}\}.
$$
Let $S^d\subset \mathbb R^{d+1}$  denote the $d$-dimensional unit sphere.
\bigskip

\medskip
\noindent{\bf Problem 2:} {\it Prove (or disprove) that ${\mathbb E}_2(A)\lesssim_{\epsilon} (\#A)^{2+\epsilon}$ if $A\subset S^2$.}
\medskip

This is known for subsets of the two dimensional paraboloid $A\subset P^2$, see \cite{bd}.
\bigskip

\medskip
\noindent{\bf Problem 3:} {\it Prove (or disprove) that ${\mathbb E}_3(A)\lesssim_{\epsilon} (\#A)^{3+\epsilon}$ if $A\subset S^1$ or $A\subset P^1$.} 
\medskip

For $S^1$, this  would follow from the unit distance conjecture. The best known unconditional bound     
is $O(\#A^{7/2})$ via the Szemer\'edi--Trotter theorem (see \cite{BombieriBourgain}).
\bigskip

In the case where all points in $A$ are 
$(\#A)^{-O(1)}$-separated all  conclusions follow from decoupling theorems in \cite{bd}.

\

\bigskip

\address {C. Demeter,
Department of Mathematics,
 Indiana University, 
Bloomington, IN 47405, USA.}

{\it E-mail address:} \email{\tt demeterc@indiana.edu}

\newpage

 \section{Functions whose Fourier transform vanishes on the sphere}

 \begin{center}
 {\it   Contributed by  Michael Goldberg  }
 \end{center}
 
\bigskip
 
\noindent{\bf Problem 4:}  {\it Let $n \geq 2$.  Does there exist a function 
$f \in L^{\frac{2n+2}{n+3}}({\mathbb R}^n)$ such that
$$\widehat{f}|_{S^{n-1}} \equiv 0$$ and 
$$\big||1-|\xi|^2\big|^{-1/2}\widehat{f}(\xi) \not \in L^2(\mathbb R^n)?$$}

\bigskip

The Stein-Tomas theorem, and approximation by Schwartz functions, implies
that restrictions of $\hat{f}(\xi)$ to $r S^{n-1}$ are continuous in $L^2(S^{n-1})$
with respect to $r>0$.  Local integrability then dictates that vanishing of $\hat{f}$ on the unit sphere
is necessary for $|1- |\xi|^2|^{-1/2}\hat{f}$ to belong to $L^2(\mathbb R^n)$.
We are curious whether it is a sufficient condition as well.  The problem is phrased to 
suggest that
perhaps it is not, and asks for a counterexample if one exists.

The more general question of finding sufficient conditions for $|1-|\xi|^2|^{-\alpha}\hat{f}
\in L^2(\mathbb R^n)$ is resolved on both sides of $\alpha = \frac12$.
For $0 \leq \alpha <\frac12$ it suffices for
$f \in L^{\frac{2n+2}{n+1+4\alpha}}(\mathbb R^n)$  \cite{BMO95}.
For $ \frac12 < \alpha < \frac32$ it also suffices for 
$f \in L^{\frac{2n+2}{n+1+4\alpha}}(\mathbb R^n)$
so long as $\hat{f} |_{S^{n-1}} \equiv 0$   \cite{G15}.

There is at least a superficial resemblance to the role of Hardy spaces in
fractional integration.  When $n=1$ the resemblance is explicit.
The sufficient condition for $|1-|\xi|^2|^{-1/2}\hat{f}$ to belong to $L^2(\mathbb R)$
is that $e^{i x}f(x), e^{-ix}f(x) \in H^1(\mathbb R)$.  This is a strictly stronger
condition than $\hat{f}(\pm 1) = 0$, and
one can construct solutions to the Problem in $L^1(\mathbb R)$ 
based on functions such as
 $$g_N(x) = \eta(x) - 2N^{-1}\cos(x)\eta(N^{-1}x)$$
and
$$h_N(x) = \eta(x-2\pi N) - \eta(x+2\pi N),$$
where $\eta$ is a Schwartz function for which
$\widehat{\eta}(\xi) = 1$ for all $|\xi| \leq 2$.

These examples don't seem to adapt well to higher dimensions, where the curvature
of $S^{n-1}$ turns into an impediment rather than an asset.

\bigskip

\bigskip

\address{M. Goldberg,
 Department of Mathematical Sciences,
University of Cincinnati,
Cincinnati, OH 45221, USA. }

{\it E-mail address:} \email{\tt Michael.Goldberg@uc.edu} 
\newpage

\section{Regularity or blow up for a 1D model of the 3D Euler equation}
\begin{center}{\it Contributed by Alexander Kiselev}\end{center}

\bigskip

The 3D Euler equation of fluid mechanics in vorticity form is given by
\begin{equation}\label{3dvort}
\partial_t \omega + (u \cdot \nabla) \omega = (\omega \cdot \nabla) u, \,\,\,u = {\rm curl} (-\Delta)^{-1}\omega, \,\,\,\omega(x,0)=\omega_0(x), \,\,\,\nabla \cdot \omega_0 =0.
\end{equation}
The equation can be set on ${\mathbb R}^3$ with decaying initial conditions, or on ${\mathbb T}^3,$ or on a domain $\Omega \subset {\mathbb R}^3,$ in which case boundary conditions
are imposed on $\partial \Omega.$
The question whether the solutions to 3D Euler equation remain globally smooth is a major open problem. One can consult \cite{MB} or \cite{MP}, for example, for
history and more details.
For the 2D Euler equation, the term on the right hand side of \eqref{3dvort} $-$ often called the vortex stretching term $-$ vanishes.
Due to this simpler structure, the global regularity for 2D Euler equation has been known since the 1930s.

Here is a very natural 1D model of \eqref{3dvort}, proposed by De Gregorio \cite{DG1,DG2}.
\begin{equation}\label{Degr}
\partial_t \omega + u \partial_x \omega = \omega \partial_x u, \,\,\, u_x = H\omega, \,\,\,\omega(x,0)=\omega_0(x).
\end{equation}
Here $H\omega = \frac{1}{\pi} p.v. \int_{{\mathbb R}} \frac{\omega(y)}{{\cdot-y}}\,\d y$ is the Hilbert transform of $\omega.$ For simplicity, let us think
of \eqref{Degr} with periodic initial data. Note that the Biot--Savart law $u_x = H\omega$ is the direct and the only natural 1D analog of $u = {\rm curl}(-\Delta)^{-1}\omega.$

Before talking about \eqref{Degr}, let us consider two simpler cases. 

\ \

{\it Case 1.}  \[ \partial_t \omega + u \partial_x \omega = 0, \,\,\, u_x = H\omega, \,\,\,\omega(x,0)=\omega_0(x). \]
This is the analog of the 2D Euler equation. Global regularity of solutions can be established similarly to that case. \\

{\it Case 2.} \begin{equation}\label{colama} \partial_t \omega  = \omega \partial_x u, \,\,\, u_x = H\omega, \,\,\,\omega(x,0)=\omega_0(x). \end{equation}
This equation has been considered by Constantin, Lax and Majda \cite{CLM}. They proved that finite time blow up can happen for certain kinds of initial data.
Moreover, the equation is explicitly solvable; by using the properties of the Hilbert transform, one can show that
\[ \omega(x,t) = \frac{4 \omega_0(x)}{(2-t H\omega_0(x))^2+ t^2 \omega_0^2(x)}.\] 
The specific properties of the Hilbert transform that are used are $H^2=-I$ and $H(\omega H \omega) = \frac12 ( (H\omega)^2 - \omega^2 ).$

\bigskip

Here are two problems connected with this theme.

\medskip\noindent
{\bf Problem 5:}  {\it 
 What happens in the full De Gregorio model \eqref{Degr}? Can the nonlinear advection term $u \partial_x \omega$ somehow balance and cancel the effect of the term that produces blow up in the
Constantin--Lax--Majda model? The mechanism for such cancellation is not clear, and so a reasonable first guess at the properties of the equation may be the possibility of finite time blow up. But in fact numerical simulations of
\eqref{Degr} suggest global regularity \cite{OSW,Sverak}. }

\newpage

\medskip\noindent
{\bf Problem 6:}  {\it  Consider an analog of the Constantin--Lax--Majda model on ${\mathbb R}^2:$
\[ \partial_t \omega = \omega {\mathcal C} \omega, \,\,\,\omega(x,0)=\omega_0(x), \]
where ${\mathcal C}$ is some Calder\'on-Zygmund operator in two dimensions. The magical properties of the Hilbert transform that allow \eqref{colama} to be solved explicitly are gone.
What can be said now? Is there still finite time blow up, and if so, how to understand its mechanism in a more robust way?
This problem (or rather, a closely related one) is stated by Constantin and Sun \cite{ConSun} in connection with some complex fluids models. }
\medskip

\bigskip

\bigskip

\address{A. Kiselev, Department of Mathematics, Rice University, 
Houston, TX 77005, USA.}

{\it E-mail address:} \email{\tt kiselev@rice.edu}

\newpage

\section{A Multi-Frequency Carleson Problem}

\centerline{\it Contributed by Ben  Krause}

\bigskip

Let $\{ \theta_1,\dots,\theta_K\} \subset \mathbb{R}$ be a collection of one-separated frequencies,
\[ \min_{ i \neq j} | \theta_i - \theta_j | \geq 1.\]

Let $P(t) := 1_{t < 0}$ and $ 1_{[-c,c]} \leq \phi \leq 1_{[-2c,2c]}, \ 1_{[-5c,5c]} \leq \varphi \leq 1_{[-10c,10c]}$ where
 $\phi$, $\varphi$ are Schwartz functions and  $0 < c \ll 1$ is sufficiently small.
Set
\[ \widehat{S_\lambda f}(\xi) := P(\xi - \lambda) \phi(\xi - \lambda)\widehat{f}(\xi).\]

The following  question  arose in work on a discrete Carleson operator along the primes \cite{CHKLP}.

\bigskip
\noindent {\bf Problem 7:} {\it What is the correct operator-norm growth in $K$ of the ``multi-frequency Carleson operator''
\[ \mathcal{C}f(x):= \sup_{|\lambda| \leq c} \left| \sum_{i=1}^K e(\theta_i x) S_\lambda ( \varphi(\xi) \widehat{f}(\xi + \theta_i) )^{\vee} (x) \right|? \]
}

\noindent Here, $e(t) := e^{2\pi i t}$. For applications, one may assume that $\{ \theta_1,\dots, \theta_K\}$ are rational frequencies each of whose denominators are comparable in size to $K$. 
\medskip

{\it Remark.}
In \cite[\S3]{CHKLP} it is shown that the $L^2$-operator norm of $\mathcal{C}$ grows like $\log^2 K$. This was accomplished by adapting a subtle argument due to Bourgain \cite{B1}, using crucially the variational Carleson theorem due to Oberlin, Seeger, Tao, Thiele, and Wright \cite{OSTTW}.

Unfortunately, away from $L^2$, the argument of Bourgain breaks down, and the most efficient (cheap) way we found to estimate $\mathcal{C}$ in $L^p, \ p \neq 2$, was through Cauchy--Schwarz, vector-valued Carleson \cite{HL}, and Rubio de Francia \cite{R}, keeping track of the loss of orthogonality near $L^1$ (see \cite{K}). This lead to the following estimates
\[ \| \mathcal{C} \|_{L^p \to L^p} \lesssim 
\begin{cases} K^{1/p} &\mbox{if } 1 < p \leq 2 \\ 
K^{1/2} & \mbox{if } 2 \leq p < \infty \end{cases}.\]

These estimates were strong enough to prove, for instance, the following discrete Carleson theorem along the primes:
\[ \left\| \sup_{k \geq 0} \left| \sum_{p \in \pm \mathbb{P}} f(n-p) \frac{e( 2^{-k} p)}{p} \log |p|  \right| \right\|_{\ell^p({\Bbb Z})} \lesssim \| f\|_{\ell^p({\Bbb Z})}, \ 3/2 < p < 4,\]
where $\pm \mathbb{P} := \{ \pm p : p \text{ is a prime} \}$. (The $\log|p|$ factor is inserted for normalization purposes.) The expected range $1 < p < \infty$ would follow from an estimate of the form
\[ \| \mathcal{C} \|_{L^p \to L^p} \lesssim_\epsilon K^{\epsilon} \]
for every $\epsilon > 0$.

\bigskip

\address{ B. Krause, Department of Mathematics, The University of British Columbia,
Vancouver, B.C.,  Canada V6T 1Z2.}

{\it E-mail address:} \email{\tt benkrause@math.ubc.ca}

\newpage

\section{Removal of Condition (R) for Fourier Restriction to Hypersurfaces in three dimensions}
\centerline{\it Contributed by Detlef M\"uller }

\bigskip

Let $S$ be a smooth hypersurface in ${\Bbb R}^3$ with Riemannian
surface measure $\d\sigma.$  Assume that $S$ is of {\it finite type},   i.e., that every tangent plane has finite order of contact with $S,$ and consider  the compactly supported measure 
$\d\mu:=\rho \d\sigma$ on $S,$ where $0\le\rho\in C_0^\infty(S).$
In my joint monograph with I. Ikromov~\cite{Muller1}, we have  investigated the possibility of Stein--Tomas type  restriction estimates of the form 
\begin{equation}\label{rest1}
\Big(\int_S |\widehat f|^2\,\d\mu\Big)^{1/2}\le C_p\|f\|_{L^p({\Bbb R}^3)},\qquad  f\in\mathcal S({\Bbb R}^3),
\end{equation}
in an almost complete way. More precisely, we have determined the optimal range of $p$'s for which  \eqref{rest1} holds true,  in terms of Newton polyhedra associated to $S,$  provided $S$ satisfies a certain condition (R). Every real-analytic surface $S$ does satisfies (R), but we do not  know if our description of this range  remains  valid  for arbitrary finite type hypersurfaces without condition (R).
\medskip

The following is a good toy model for testing this question.

\noindent{\it Example.}\label{ex1+}
{\rm Let $S$ be the graph of  the function 
 $$
 \phi (x_1,x_2):=(x_2-x_1^{m})^n, \qquad n,m\ge 2,
 $$
over a small neighborhood $U$ of the origin in ${\Bbb R}^2.$ 
The coordinates $(x_1,x_2)$ are not adapted to the function $\phi$  in the sense of Varchenko,  and the so-called Newton distance  of $\phi$ in these coordinates $(x_1,x_2)$ is given by 
$$ d:=\frac {nm}{m+1}.$$ 
Adapted coordinates are 
$y_1:=x_1, y_2:=x_2-x_1^m,$ in which $\phi$ is given by
 $\phi^a(y_1,y_2)=y_2^n.$
As explained in our papers,  from the (extended) Newton polyhedron of ${\phi^a}$ one can compute the so-called $r$-height $h^r$  of $\phi,$ and our main theorem implies that \eqref{rest1} holds true if and only if $p'\ge p'_c=2(h^r+1).$ In this example, we  find that  the critical conjugate exponent $p'_c$ is given by
$$ p'_c=\left\{  \begin{array}{cc}
2d+2, & \mbox{   if } n\le m+1 ,\hfill\\ 
2n, & \mbox{   if } n>m+1.\hfill 
\end{array}\right.
$$

\medskip
{\bf Problem 8:} {\it Does this result remain valid if we perturb $\phi$ by a smooth function which is flat at the origin (provided we choose $U$ sufficiently small), for instance for a perturbed function 
$$
\phi_1(x_1,x_2):=(x_2-x_1^{m})^n +{\varphi}(x_1),
$$
where  all derivatives of ${\varphi}$ vanish at the origin? For non-trivial ${\varphi},$ such a function $\phi_1$ does not satisfy our condition (R), in contrast with $ \phi_2(x_1,x_2):=(x_2-x_1^{m}-{\varphi}(x_1))^n.$}
\medskip

{For more information on this question the reader is referred to 
\cite{Muller1}, \cite{Muller2}, \cite{Muller3}.}

\bigskip
\address{D. M\"uller,
Mathematisches Seminar,
Christian-Albrechts-Universit\"at Kiel,
D-24118 Kiel, Germany.} 

{\it E-mail address:} \email{{\tt mueller@math.uni-kiel.de}}

 \newpage

\section{Asymptotic behavior of the $N$-th polarization constant for ``fat'' Cantor sets}
\begin{center}{\it Contributed by Alexander Reznikov}\end{center}

\bigskip

Let $\mathcal{C}\subset [0,1]$ denote a ``fat'' Cantor set, where ``fat'' means $m_1(\mathcal{C})>0$, and $m_1$ denotes the Lebesgue measure on $[0,1]$. For the sake of this problem, one can take any particular Cantor set. Fix a number $s>1$.

By $\omega_N$ we denote a set of points with possible repetitions, with cardinality $\#\omega_N = N$, counting multiplicities. For any such $\omega_N$ denote
$$
P_s(\mathcal{C}; \omega_N):=\inf_{y\in \mathcal{C}} \sum\limits_{x_j\in \omega_N} \frac{1}{|y-x_j|^s}.
$$
Finally, denote
$$
\mathcal{P}_s(\mathcal{C}; N):=\sup_{\omega_N\subset \mathcal{C}} P_s(\mathcal{C}; \omega_N).
$$

It is well known, see \cite{ES}, that there exist positive constants $c_1$ and $c_2$, such that
$$
c_1N^s \leqslant \mathcal{P}_s(\mathcal{C}; N)\leqslant c_2 N^s, \;\;\; N\geqslant 2.
$$

\medskip\noindent
{\bf Problem 9:}  {\it Prove that 
$$
\lim_{N\to \infty} \frac{\mathcal{P}_s(\mathcal{C}; N)}{N^s}
$$
exists. Further, prove (or disprove) 
$$
\lim_{N\to \infty} \frac{\mathcal{P}_s(\mathcal{C}; N)}{N^s} = \frac{2(2^s-1)\zeta(s)}{m_1(\mathcal{C})^s},
$$
where $\zeta(s)$ is the classical Riemann zeta-function. 
}
\medskip

We remark that for the circle $\mathbb{T}$ instead of $\mathcal{C}$, it is known that the limit is equal to
$$
\lim_{N\to \infty} \frac{\mathcal{P}_s(\mathbb{T}; N)}{N^s} = \frac{2(2^s-1)\zeta(s)}{(2\pi)^s}.
$$
We believe that the numerator should remain the same, and that the denominator is always the ``length'' of our set, raised to the power $s$.

\bigskip

\bigskip

\address{A. Reznikov, 
Department of Mathematics, 
Vanderbilt University, 
Nashville, TN 37240, USA}

{\it E-mail address:} \email{\tt aleksandr.b.reznikov@vanderbilt.edu}

\newpage

\section{Three problems in Harmonic Analysis}

\centerline{\it Contributed by Andreas Seeger}

\bigskip

\bigskip

\subsection{Convolutions with radial kernels}

Let $\sigma_r$ be the surface measure on the sphere of radius $r$ in $\mathbb R^d$, $d\ge 2$, and let $\eta$ be a Schwartz function in $\mathcal S(\mathbb R^d)$ 
whose Fourier transform is zero on a neighborhood of the origin.

 \medskip

\noindent
{\bf Problem 10:} {\it For $1<p<2$, $d> \frac{p}{2-p}$, 
\begin{equation}\label{conj1}
\Big\|\int_1^\infty (\eta*\sigma_r * g)(\cdot, r) \d r \Big\|_{L^p(\mathbb R^d)}\le C(p,d,\eta) \Big(\int_1^\infty\int_{\mathbb R^d} \big|g(x,r)|^p \, \d x\, r^{d-1} 
\d r\Big)^{\frac1p}.
\end{equation}}

If one applies  the inequality to functions of the form
$g(x,r)=f(x) k(r)$, one can derive effective characterizations 
of  radial multipliers 
of $\mathcal F L^p(\mathbb R^d)$.  In particular, if $K$ is radial and $\widehat K$ is compactly supported away from the origin, then
$\|\widehat K\|_{M^p}\approx \|K\|_{L^p}$. More generally, let $\phi$ be a nontrivial $C^\infty_c(\mathbb R^+)$ function. Then \eqref{conj1} implies 
$$\|m(|\cdot|)\|_{M^p(\widehat {\mathbb R^d})} \approx  \sup_{t>0}\big \|\mathcal F^{-1}[\phi(|\cdot|)  m(t|\cdot|)] \big\|_{L^p(\mathbb R^d)}.$$
Thus,  for such $p$, a convolution operator with radial kernel is bounded on $L^p(\mathbb R^d)$ if and only if  it is bounded on the subspace $L^p_{\text{rad}}(\mathbb R^d)$ consisting of radial $L^p$ functions, cf. \cite{gs}.
These  statements are  false if $2\le d\le \frac p{2-p} $ 
(as then $\mathcal F L^p(\mathbb R^d)$ contains unbounded radial  functions which are compactly supported away from the origin). 

Inequality \eqref{conj1} 
is also  known to imply an  endpoint version of the so-called local smoothing conjecture (see \cite{sogge} for the wave equation, 
and for an endpoint estimate for cone multipliers in $\mathbb R^{d+1}$, see \cite{hns1},  \cite{hns2}).
Consult  \cite{taorestr}, \cite{taowt}  for   implications on the 
Bochner--Riesz, Fourier restriction,  Kakeya--Nikodym,  and other related  conjectures.

Note that  \eqref{conj1} is trivial for $p=1$. Given  $1<p<2$,  it was shown in \cite{hns1} that  \eqref{conj1} holds in  dimensions
$d>\frac{p+2}{2-p}$. This gives a known nontrivial range of $p$ in dimensions $d>3$. See also \cite{cladek} on improvements of \eqref{conj1} for special  classes of functions.
Weaker forms of the conjecture are known to hold in larger ranges,
for example 
\begin{equation}\label{conj2}
\Big\|\int_1^\infty (\eta*\sigma_r * g)(\cdot, r) \d r \Big\|_{L^p(\mathbb R^d)} 
\le C_\varepsilon(p,d,\eta) \Big(\int_1^\infty\int_{\mathbb R^d} \big|g(r,x)|^p \, \d x\, r^{d-1+\varepsilon} 
\d r\Big)^{\frac1p}
\end{equation}
holds for all $\varepsilon>0$ if $d\ge \frac{3p-2}{2-p}$ (of course for fixed $\varepsilon$ the range is  better, depending on $\varepsilon$). This follows from the validity of the Wolff $\ell^{p'}$-decoupling inequalities (\cite{wolff}, \cite{LW})  in the optimal   range  $p'>\frac{2(d+1)}{d-1}$ which is a consequence of the more recent results  by Bourgain and Demeter \cite{bd}. In particular 
\eqref{conj2} holds for all $\varepsilon >0$ for the  nontrivial ranges $p\in [1,\frac 65]$ in dimension two,  and $p\in [1, \frac 43]$ in dimension three.

  \subsection{Hyperbolic cross maximal function}
  Let $\Phi$ be a $C^\infty_c$  function on the real line and let $\Phi(0)\neq 0$.
 Define the operator $T_t$ acting on functions in $\mathbb R^2$ by
 $\widehat {T_t f} (\xi) = \Phi( t^{-1} \xi_1^2\xi_2^2) \widehat f(\xi) $. The operators $T_t$ are bounded on $L^p(\mathbb R^2)$, with the same operator norm.
 
 \smallskip
 \noindent
 {\bf Problem 11:} {\it Does the hyperbolic cross maximal inequality 
 \begin{equation}\label{hyp}\big \|\sup_{k\in \mathbb Z}  |T_{2^k} f |\big\|_{L^p(\mathbb R^2)} \le C\|f\|_{L^p(\mathbb R^2)}
 \end{equation} hold for some $p$?}
 
  \smallskip
 
 If \eqref{hyp} held for some $p_0\in (1,\infty)$, then \eqref{hyp} would hold for all $p\in (1,\infty)$, by standard Calder\'on--Zygmund techniques, see e.g. \cite{dt}. 
 An affirmative answer would imply an affirmative answer to a question posed by Dappa and Trebels  in \cite{dt}, concerning  the hyperbolic  Riesz means $H^\lambda_s $ defined by $$\widehat { H^\lambda_s f }(\xi)=\Big(1-\frac{\xi_1^2\xi_2^2}{s^2}\Big)^\lambda_+\widehat f(\xi).$$  For sufficiently large $\lambda$, and all $f\in L^p(\mathbb R^2)$,   does 
 $H^\lambda_s f(x) $ converge to $f(x)$ almost everywhere, as $s\to \infty$?

  \subsection{Rough singular integrals on $L^1$}
 Let $d\ge 2$, $q>1$, $\Omega\in L^q(S^{d-1})$,  and assume $\int_{S^{d-1} }\Omega(x')\d\sigma(x')=0$.
 Let $0<\varepsilon<1$ and define the family of rough singular integral operators $\{T_\varepsilon\}$   by 
 $$T_\varepsilon f(x) = \int_{\varepsilon<|y|\le 1}  |y|^{-d} \Omega( \tfrac y{|y|}) f(x-y) \d y.$$
 It was shown in \smallskip\cite{seeger-rough} that the operators $T_{\epsilon}$ are weak type $(1,1)$ with constants independent of $\varepsilon$, see also \cite{christ-rdf}, \cite{hofmann}  for earlier  results in low dimensions, 
and \cite{taorough} for an extension to stratified groups.
 
 \smallskip
  \noindent
 {\bf Problem 12:} {\it  For  $f\in L^1(\mathbb R^d)$,  does $\lim_{\varepsilon\to 0} T_\varepsilon f(x)$ exist for almost every $x\in \mathbb R^d$?

Equivalently, does
$\sup_{0<\varepsilon<1} |T_\varepsilon f(x)|$ define an operator of weak type $(1,1)$?
}

\

This problem has been around for a while and was proposed again in  \cite{seeger-cj}. Andrei Lerner (personal communication, El Escorial 2016) pointed out renewed interest in it,  in connection with problems on control of rough singular integral operators via  sparse operators, and related questions about weighted norm inequalities, cf. \cite{lerner}.

\bigskip

\bigskip

\address{A. Seeger,  Department of Mathematics,  University of Wisconsin-Madison,   Madison, WI 53706, USA.}

{\it E-mail address:} \email{\tt seeger@math.wisc.edu}

\newpage

\section{Triangular Hilbert transform}

\centerline{\it  Contributed by  Christoph Thiele }

\bigskip

\maketitle

The triangular Hilbert transform is given as the principal
value integral 
$$\Lambda_3(f,g,h)=p.v. \iiint f(x,y)g(y,z)h(z,x)\frac 1{ x+y+z}\, \d x\d y\d z,$$
where $f,g,h$ are three test functions in the plane.

\medskip
\noindent{\bf Problem 13:}  {\it 
Prove the a priori estimate
$$|\Lambda_3(f,g,h)|\le C \|f\|_{L^3}\|g\|_{L^3}\|h\|_{L^3},$$
with a universal constant $C$. }
\bigskip

The problem will be solved if one proves any similar 
estimate with Lebesgue exponents $3$ replaced by $p,q,r$ on the right hand side,
which then by homogeneity of the problem have to satisfy $1/p+1/q+1/r=1$.
Namely, one will have the symmetric estimate under permutations of the exponents, and interpolation will yield the central estimate with all exponents $3$.

The triangular Hilbert transform is one of a family of simplex Hilbert 
transforms, one for each degree of multilinearity. The adjacent degrees are
$$\Lambda_2(f,g)=p.v. \iint f(x)g(y)\frac 1{ x+y}\, \d x\d y ,$$
\begin{align*}
& \Lambda_4(f,g,h,k) \\
& =p.v. \iiiint f(x,y,z)g(y,z,u)h(z,u,x)k(u,x,y)\frac 1{ x+y+z+u}\, \d x\d y\d z\d u .
\end{align*}
It turns out $\Lambda_2$ is equivalent to the classical Hilbert transform with the well known and relatively easy estimate 
$$|\Lambda_2(f,g)|\le C \|f\|_{L^2}\|g\|_{L^2},$$
while the estimate 
$$|\Lambda_4(f,g,h,k)|\le C \|f\|_{L^4}\|g\|_{L^4}\|h\|_{L^4}\|k\|_{L^4}$$
appears well out of reach of present technology.
The case $\Lambda_3$ also requires some new ideas, but at least there is an
encouraging variety of recent results in the vicinity of the problem which
seem to suggest some helpful technology and possible routes of attack.

There is a dyadic model for the triangular Hilbert transform.
A dyadic interval is an interval of the form $I=[2^kn, 2^k(n+1))$  
with integers $k,n$ where we assume $n$ non-negative. 
If $I_1,I_2,I_3$ are three such dyadic intervals of the same
length, we write $I_1\oplus I_2 \oplus I_3=0$ if $n_1\oplus n_2\oplus n_3=0$,
where addition is that of the group $(\mathbb{Z}/2\mathbb{Z})^{\mathbb{N}}$ identified with the natural 
numbers via binary expansion of the latter. The dyadic model is then
$$\sum_{k\in \mathbb{Z}} 2^{-k} \epsilon_k \sum_{|I|=|J|=|K|=2^k, I\oplus J\oplus K=0} \iiint f(x,y)g(y,z) h(z,x)h_I(x)h_J(y)h_K(z) \d x\d y\d z, $$
 where $\epsilon_k$ are numbers in $[-1,1]$ and $h_I$ is the Haar function
associated to the dyadic interval $I$, normalized so that $\|h_I\|_{L^\infty}=1$.
If all $\epsilon_K$ are equal to $1$, the form telescopes into the 
integral of the pointwise product $fgh$. Proving $L^3$ bounds as in
the above problem for this model operator is expected to be somewhat
easier than the main problem, but be significant in terms of 
principle ideas for the attack on the main problem. 

Of the many papers on the triangular Hilbert transform and related
problems, we single out the recent ones \cite{ktz} and \cite{dkt} and 
refer for further contributions to the chain of references therein.

 The paper \cite{ktz} proves the desired bounds for the dyadic model of the triangular Hilbert transform under crucially simplifying structural assumptions on one of the 
three functions $f,g,h$, such as for example elementary tensor structure. 
Two types of structural assumptions are discussed there, leading
to connections with objects of study in time-frequency analysis such as Carleson's operator controlling almost everywhere convergence of Fourier series and the bilinear Hilbert transform. It is expected that for the main problem time-frequency analysis has to be replaced by a more general theory that is not based on Fourier analysis but on more general expansions. Speculation goes as far as the circle of ideas around the recently solved Kadison--Singer problem \cite{mss} and its consequences for example on Riesz frames.

The paper \cite{dkt} proves bounds on the truncated triangular Hilbert transform
$$\Lambda_{3,N}(f,g,h)=\iiint _{1<|x+y+z|<N} f(x,y)g(y,z)h(z,x)\frac 1{ x+y+z}\, \d x\d y\d z ,$$
with bounds on the growth of the constant that are in between trivial and desired, namely
$$|\Lambda_{3,N}(f,g,h)|\le C (\log N)^{1/2} \|f\|_{L^3}\|g\|_{L^3}\|h\|_{L^3}.$$
The quest for such growth estimates goes back to \cite{tt} and \cite{pz}, trivial being the power
$(\log N)^1$ and desired being $(\log N)^0$.
More generally, the paper \cite{dkt} proves similar estimates for all simplex Hilbert transforms, with less improvement on the exponent as the degree of multilinearity increases. The technique in this paper is called ``twisted technology'', which also recently played a role in ergodic theory \cite{dkst}, proving quantitative bounds on convergence of ergodic means for two commuting transformations.

\bigskip

\address{C. Thiele, 
Mathematisches Institut 
der Universit\"at Bonn, Endenicher Allee 60,
 53115 Bonn, Germany.}
 
{\it E-mail address:}  \email{{\tt thiele@math.uni-bonn.de}} 
 
 \newpage
 \section{Problems  on lacunary maximal functions}
\centerline{\it Contributed by James Wright}

\bigskip

We  consider a compactly supported finite Borel measure $\mu$ on ${\mathbb R}^d$ 
and define its dyadic dilates by
$\inn{\mu_k}{f}= \inn{\mu}{f(2^k\cdot)}$. We present here three open
problems about the corresponding lacunary maximal function  given by
$${\mathscr {M}} f(x)= \sup_{k\in {\mathbb {Z}}} |f*\mu_k(x)|.$$
The  dilates  $2^k$ can be replaced by
more general lacunary dilates $\lambda_k$
satisfying $\inf_k \lambda_{k+1}/\lambda_k > 1$. For more details and other open problems, see
\cite{sw}.

\vskip 10pt

For the first two problems, suppose that 
$\mu$  satisfies the Fourier decay condition
\begin{equation}\label{decay}
\bigl| \widehat \mu(\xi) \bigr| \  \le \ C \, |\xi|^{-\varepsilon} \ \ {\rm for \ some} \ \varepsilon>0. 
\end{equation}

\vskip 15pt

{\bf Problem 14}: {\it Does 
${\mathscr {M}}$ satisify a weak-type $(1,1)$ bound? That is, does ${\mathscr {M}} : L^1({\mathbb R}^d) \to L^{1,\infty}({\mathbb R}^d)$?}
\vskip 10pt\
To the best of our knowledge no counterexample and no example
is known for the case when $\mu$ is a singular measure with the decay assumption \eqref{decay} on $\widehat \mu$.
We recall that the $L^p$-boundedness of ${\mathscr {M}}$ holds in the
range $1<p \le \infty$, see {\it e.g.} \cite{duordf}.

\vskip 15pt

{\bf Problem 15}: {\it Does ${\mathscr {M}}$ map the Hardy space $H^1({\mathbb R}^d)$ to $L^{1,\infty}({\mathbb R}^d)$? }
\vskip 10pt
This is a more accessible problem and it is known to hold for $\mu$ being surface measure on the unit
sphere ${S}^{d-1}$. This is due to M. Christ \cite{christ} and represents the first
such result for a singular measure $\mu$ satisfying \eqref{decay}. More recently Christ's result has been
extended to any $\mu$ satisfying a dimensional assumption on ${\rm supp}(\mu)$, coupled with an
optimal $L^p$ smoothing assumption near $p=1$ for the corresponding averaging operator $f \to f*\mu$.
See \cite{sw}. These assumptions are known to hold for surface measure on any hypersurface
on which the Gaussian curvature does not vanish to infinite order. They also hold 
for arclength measure on a compact curve in ${\mathbb R}^3$
with nonvanishing curvature and torsion. 

\vskip 10pt

For the third and final problem, we examine lacunary maximal functions associated to measures
$\mu$ which do not necessarily satisfy the Fourier decay condition \eqref{decay}. 

\vskip 5pt

Let $\Omega$ be a convex open domain in the plane with compact closure so that the origin is contained in $\Omega$.
We let  $\sigma$ be the arclength measure on the  boundary ${\partial\Omega}$ and consider the lacunary maximal operator associated to ${\partial\Omega}$,
 $$
{\mathcal {M}} f(x) =\sup_{k\in {\mathbb {Z}}} \Big|\int_{{\partial\Omega}}
 f(x- 2^{k} y) \d\sigma(y) \Big|.
$$

\vskip 15pt

{\bf Problem 16}: {\it For any fixed $1<p<\infty$, give a necessary and sufficient geometric condition on ${\partial\Omega}$ 
which guarantees ${\mathcal {M}}$ is bounded on 
$L^p({\mathbb R}^2)$. Below we suggest a possible geometric condition.}

\vskip 10pt

As mentioned above,  if
$|\widehat \sigma(\xi)|=O(|\xi|^{-\varepsilon})$ for some $\varepsilon>0$, then ${\mathcal {M}}$ is  $L^p$ bounded for all 
$p$ in the range $1<p\le \infty$. Now  we
are aiming for much weaker conditions.

\vskip 3pt

The decay of $\widehat \sigma$  is strongly related to a geometric
quantity. Given a unit vector $\theta$, let $\ell^+(\theta)$ be the unique
supporting line with $\theta$ an outer normal to ${\partial\Omega}$, i.e.
the affine line perpendicular to $\theta$ which intersects ${\partial\Omega}$
so that $\Omega$ is a subset of the halfspace $\{x: x=y-t\theta: t>0, y\in\ell^+(\theta)\}$. Similarly define
$\ell^-(\theta)$ as the unique
affine line perpendicular to $\theta$ which intersects ${\partial\Omega}$
and  $\Omega$ is a subset of the halfspace
$\{x: x=y+t\theta: t>0, y\in\ell^-(\theta)\}$.
For small $\delta>0$, define the arcs (or \lq caps\rq )
$${\mathcal {C}}^{\pm}(\theta,\delta) =\{y\in{\partial\Omega}: {\text{\it dist}}(y, \ell^{\pm})\le \delta\}.$$
By compactness considerations it is easy to see that there exists $\delta_0>0$ so that, for all $\theta\in S^1$ and all $\delta<\delta_0$, the  arcs
${\mathcal {C}}^{+}(\theta,\delta) $
 and ${\mathcal {C}}^{-}(\theta,\delta)$ are disjoint.
Let  $\Lambda(\theta,\delta)$ be the maximum of the length of these caps:
$$\Lambda(\theta,\delta) = \max_\pm \sigma({\mathcal {C}}^{\pm}(\theta,\delta))\,.$$
The analytic  significance of this quantity is that it gives  a
 very good estimate   for the size of Fourier transform $\widehat\sigma $,
namely,
for every $\theta\in S^1$ and $R\ge 1$,
$$
|\widehat \sigma(R\theta)| \le C_\Omega \Lambda(\theta, R^{-1}).
$$
This is shown in \cite{bnw} under the hypothesis for convex domains in the plane with smooth boundary, with no quantitative assumption on the second derivative. The general case follows by a simple approximation procedure.

In \cite{sw}, the following theorem was established, giving a solution to Problem 16 when $p=2$.

\medskip 

\noindent {\bf Theorem.} {\it The operator ${\mathcal {M}}$ is bounded on $L^2({\mathbb {R}}^2)$ if and only if
$$
\sup_{\theta\in S^1} \int_0^{\delta_0} \Lambda (\theta,\delta) ^2 \frac{\d\delta}{\delta}
\ < \ \infty.
$$} 

By testing ${\mathcal {M}}$ on functions supported in thin strips, one 
obtains the  necessary condition
$$
\sup_{\|f\|_{L^p}=1}\|{\mathcal {M}} f\|_{Lp} \ge c
\sup_{\theta\in S^1}  \Big(\int_0^{\delta_0} \Lambda (\theta,\delta) ^p \frac{\d\delta}{\delta}\Big)^{1/p}
$$
for any $p$.
We note  that if ${{\partial\Omega}}$ has only one ``flat'' point near which ${{\partial\Omega}}$
can be parametrized as the graph of $C+\exp(-1/|t|^a)$ with $C\neq 0$,
then this $L^p$ condition
holds if and only if $a<p$. Thus in this case  $L^2$ boundedness of ${\mathcal {M}}$
 holds if and only if $a<2$.

\vskip 15pt

\noindent
{\bf Problem 16 (precise form):} {\it For any $p\not= 2$, is the geometric condition
$$
\sup_{\theta\in S^1} \int_0^{\delta_0} \Lambda (\theta,\delta) ^p \frac{\d\delta}{\delta}
\ < \ \infty
$$
sufficient for ${\mathcal {M}}$ to be bounded on $L^p({\mathbb R}^2)$?}

\bigskip

\bigskip

\address{J.  Wright,
Maxwell Institute for Mathematical Sciences 
 and the School of Mathematics,  University of Edinburgh, 
 Edinburgh EH9 3FD, Scotland.}

{\it E-mail address:}  \email{\tt J.R.Wright@ed.ac.uk }

\newpage

{\it Editors:}

\end{document}